\documentstyle[12pt]{amsart}

\newcommand{\labell}[1] {\label{#1}}

\numberwithin{equation}{section}

\newtheorem{lemma}{Lemma}
\newtheorem{theorem}{Theorem}

\newtheorem{def-theo}[theorem]{Definition-Theorem}

\theoremstyle{definition}
\newtheorem*{Definition}{Definition}
\newtheorem*{Remark}{Remark}

\newenvironment{Remarks}
{\subsection*{Remarks} \begin{enumerate}}{\end{enumerate}\smallskip}

\expandafter\chardef\csname pre amssym.def at\endcsname=\the\catcode`\@ 
\catcode`\@=11 
\def\undefine#1{\let#1\undefined} 
\def\newsymbol#1#2#3#4#5{\let\next@\relax 
 \ifnum#2=\@ne\let\next@\msafam@\else 
 \ifnum#2=\tw@\let\next@\msbfam@\fi\fi 
 \mathchardef#1="#3\next@#4#5}
\def\mathhexbox@#1#2#3{\relax 
 \ifmmode\mathpalette{}{\m@th\mathchar"#1#2#3}%
 \else\leavevmode\hbox{$\m@th\mathchar"#1#2#3$}\fi} 
\def\hexnumber@#1{\ifcase#1 0\or 1\or 2\or 3\or 4\or 5\or 6\or 7\or 8\or 
 9\or A\or B\or C\or D\or E\or F\fi} 
 
\font\teneufm=eufm10 
\font\seveneufm=eufm7
\font\fiveeufm=eufm5 
\newfam\eufmfam
\textfont\eufmfam=\teneufm 
\scriptfont\eufmfam=\seveneufm 
\scriptscriptfont\eufmfam=\fiveeufm 
\def\frak#1{{\fam\eufmfam\relax#1}} 
\catcode`\@=\csname pre amssym.def at\endcsname 

\def	\red	{{\operatorname{red}}}

\def\codim{\hbox{\normalshape codim}\,}

\def	\CC	{{\Bbb C}}
\def	\RR	{{\Bbb R}}
\def	\CP	{{\Bbb C}{\Bbb P}}

\def	\LL	{\boldsymbol L}

\def	\RRoch	{\operatorname{RR}}

\begin{document}


\setlength{\smallskipamount}{6pt}
\setlength{\medskipamount}{10pt}
\setlength{\bigskipamount}{16pt}





\title[Cobordism theory and localization formulas]
{Cobordism theory and localization formulas for Hamiltonian group actions}

\author[Viktor Ginzburg]{Viktor Ginzburg}
\author[Victor Guillemin]{Victor Guillemin}
\author[Yael Karshon]{Yael Karshon}

\address{V. Ginzburg: Department of Mathematics, UC Berkeley,
Berkeley, CA 94720}
\email{ginzburg@@math.berkeley.edu}

\address{V. Guillemin: Department of Mathematics, MIT, Cambridge, 
MA 02139-4307}
\email{vwg@@math.mit.edu}

\address{Y. Karshon: The Hebrew University of Jerusalem, Institute of
Mathematics, Givat Ram, 91904 Jerusalem, Israel}
\email{karshon@@math.huji.ac.il}

\thanks{Available electronically from dg-ga/9601003}
\date{December 1995}

\thanks{The work of all three authors was partially supported by the NSF}

\maketitle

\bigskip

In this note we will report on some symplecto-geometric\footnote{
However, we will usually allow our closed two-forms to be
degenerate.
}
applications of the following result: 

\begin{theorem} \labell{theorem1}
Let $M$ be a compact oriented $2d$-dimensional manifold on which the group
$S^1$ acts. Suppose that this action is quasi-free and has finitely many
fixed points. Then $M$ is cobordant to a disjoint union of $N$
copies of $\CP^d$, where $N$ is the number of fixed points.
\end{theorem}

A more detailed discussion of this theorem and its applications
will appear in  \cite{GGK}.\footnote{A similar cobordism technique has 
been developed and used independently of us
by Shaun Martin in his study of Hamiltonian group actions on symplectic
manifolds \cite{Mar}. We wish to thank him for acquainting us with 
his results in numerous fruitful discussions which have been of a great
value to us.}

\begin{Remark}
``Quasi-free'' means that $S^1$ acts freely on the
complement of the fixed point set.
The cobordism is an $S^1$-equivariant cobordism
with the $\CP^d$'s carrying projective $S^1$-actions. Furthermore,
the cobordism is oriented although the $\CP^d$'s need not be given their
standard orientations. 
\end{Remark}

\begin{pf} 
The circle $S^1$ acts on $M \times \CC$ by the product of 
its action on $M$ and its standard action on $\CC$. 
The fixed points of this action are
$$
q_k = (p_k, 0)\ ,\ p_k \in M^{S^1}\ ,
$$
where $M^{S^1}$ is the set of fixed points for the circle action on $M$. 
Denote by $U_k$ an $S^1$-invariant open ball around $q_k$ (with respect to
some invariant Riemann metric).  Let $W$ be the subset of $M \times \CC$
obtained by excising the $U_k$'s and the set $|z| >1$. Since $S^1$ acts
freely on $W$, the quotient $W/S^1$ is a compact manifold-with-boundary,
and, modulo orientations,
\begin{equation}\labell{1}
\partial (W/S^1) = M \cup \coprod_k (\partial U_k)/S^1\ .
\end{equation}
Let $T_{q_k}$ be the tangent space of $M \times \CC$ at $q_k$.
The linear isotropy action of $S^1$ on $T_{q_k}$ is free except at the
origin, hence there is an $\RR$-linear identification
\begin{equation}\labell{2}
T_{q_k} \cong \CC^{d+1}
\end{equation}
which converts this action into the action ``multiplication by
$e^{i\theta}$''. 
Via \eqref{2} one can identify $U_k$ with the set $||z|| < \varepsilon$ 
and hence identify $\partial U_k/S^1$ with $\CP^d$. Thus, by
\eqref{1}, $M$ is cobordant to the disjoint union of $N$ copies of 
$\CP^d$. 

If the isomorphism \eqref{2} respects orientation, the
$k$-th $\CP^d$ is equipped with the complex orientation;
otherwise we take the opposite orientation.

The action of $S^1$ on the first component of $M\times \CC$ commutes with 
the diagonal action and thus descends to $W/S^1$, making it into an
equivariant cobordism.
\hfill\end{pf}

In what follows we will exploit the fact that certain structures on $M$
which are preserved by the $S^1$-action (for instance, a closed 2-form,
a stable complex structure, or a $G$-action) can be incorporated into
the cobordism \eqref{1}. For example, we have already seen how the cobordism
inherits the $S^1$-action.
Before introducing these new structures, though, we would like to
remark that the assumptions made in Theorem \ref{theorem1} are
unrealistically strong. What happens if the action of $S^1$ on $M$
is not quasi-free, or does not have finitely many fixed points?
If the action is not quasi-free, the induced action on $W$ is 
locally free, but not free, so the cobording space $W/S^1$ is now 
an {\it orbifold}. Moreover, under the 
identification \eqref{2}, the $S^1$-action on $\CC^{d+1}$ has the form
\begin{equation}\labell{3}
	e^{i\theta}z = 
	(e^{i m_1 \theta} z_1, \ldots, e^{i m_{d+1} \theta} z_{d+1}).
\end{equation}
Here the $m_j$'s are nonzero
integers which are no longer equal to $\pm 1$, and the
quotient, $\partial U_k / S^1$, is also an orbifold (a twisted projective
space). Therefore, from \eqref{1} one gets an equivariant 
{\it orbifold} cobordism between $M$ and a disjoint union of twisted 
projective spaces. If, on the other hand, the action is quasi-free 
but the fixed point set is not finite, one can prove:

\begin{theorem} Let $X_k$, $k=1, \ldots N$, be the connected components 
of the fixed point set $M^{S^1}$. Then there is an equivariant cobordism
\begin{equation}\labell{4}
M \sim \coprod_{k=1}^N B_k
\end{equation}
and fibrations
\begin{equation}\labell{5}
\CP^{m_k} \hookrightarrow B_k \to X_k
\end{equation}
where $2 m_k = \codim X_k$.
\end{theorem}

Finally, if neither of these hypotheses holds, there exists an
orbifold cobordism of the form \eqref{4}; however, the $\CP^{m_k}$'s
in \eqref{5} have to be replaced by twisted projective spaces.

The main result of this article is a reformulation of Theorem
\ref{theorem1} in the setting of Hamiltonian group actions.
To state this result we will first explain what is
meant by ``cobordism'' in this setting. 

\begin{Definition}
Let $G$ be a compact Lie group and let $(M_r, \omega_r)$ be a
compact oriented\footnote{ 
N.B. The orientation of $M$ need {\it not} be identical with its 
symplectic orientation (see [Gin]).
} 
$2d$-dimensional symplectic manifold on which $G$ acts in a Hamiltonian
fashion, with $\phi_r \colon M_r \to \frak g^*$ as the associated moment
mapping. One says that $(M_1, \omega_1, \phi_1)$ and 
$(M_2, \omega_2, \phi_2)$ are cobordant as Hamiltonian $G$-spaces if there
exists a compact oriented $(2d+1)$-dimensional manifold-with-boundary
$W$, a closed two-form $\omega$, and a Hamiltonian action of $G$ with
moment map $\phi \colon W \to \frak g^*$, such that
$$ \partial W = M_1 \cup (-M_2) $$ 
and such that the pull-backs of $\omega$ and $\phi$ to $M_r$ are $\omega_r$ 
and $\phi_r$.
\end{Definition}

\begin{Remarks}
\item
This definition is different from the definition used in \cite{Gin}
in that we do {\em not\/} impose any rank or non-degeneracy condition
on the closed form $\omega$. 

\item
By a ``Hamiltonian space" we will generally mean a manifold
equipped with a group action, a closed invariant two-form
(not necessarily symplectic), and a corresponding moment map.
Our definition of cobordism extends word-for-word to the case 
where the two-forms $\omega_r$ are degenerate, and therefore
we have the notion of cobordism of Hamiltonian spaces.

\item
The equivariant form $\omega_r+\phi_r$ can be altered by an equivariant 
coboundary without changing the cobordism class of the Hamiltonian
space $(M_r, \omega_r, \phi_r)$, so the cobordism class of a 
Hamiltonian $G$-space depends only on the underlying $G$-manifold
and on the cohomology class $[\omega_r + \phi_r]$. 

\item
In all our applications the group $G$ will be a torus.
\end{Remarks}

Apropos of this definition we have proved \cite{GGK} the following:

\begin{lemma}
\labell{lemma:regular}
Let $T$ be a torus and let $(W,\omega,\phi)$ be a cobordism
of Hamiltonian $T$-spaces between $(M_r,\omega_r,\phi_r)$, $r=1,2$.
Assume that $a$ is a regular value of $\phi_r$ for $r=1, 2$. 
Then one can perturb $(\omega, \phi)$ 
away from the boundary of $W$
so that $a$ becomes a regular value of $\phi$.
\end{lemma}

A consequence of this is that {\em cobordism commutes with reduction}.
More explicitly, let $Z = \phi^{-1}(a)$. 
The fact that $a$ is a regular value of $\phi$ implies that
the action of $T$ on $Z$ is locally free, hence $Z/T$ is an
orbifold-with-boundary. There exists a closed two-form $\omega^\red$ 
on $Z/T$ whose pull-back to $Z$ is equal to the restriction
to $Z$ of $\omega$. The boundary components of $Z/T$ are the reduced 
spaces $M_r^\red$, and the restrictions of $\omega^\red$ to each of 
these boundary components is the usual reduced symplectic form. 
All of this remains true even when $\omega$ is degenerate, and implies:

\begin{theorem} \labell{theorem4}
Let $a$ be a regular value of $\phi_r$ for $r=1,2$. Then
the reduced spaces $M_r^\red = \phi_r^{-1}(a)/T$ are cobordant 
(as symplectic orbifolds).
\end{theorem}

This theorem indicates that, when working with Hamiltonian group actions, 
orbifold cobordism is a far more natural notion than the usual notion of 
cobordism.  From now on we will tacitly allow all our cobordisms to be 
orbifold cobordisms.
We would like to note, however, that not every orbifold is cobordant 
to a manifold (even over the rationals), and in passing to orbifolds 
we essentially enlarge the set of cobordism classes \cite{GGK}.

Before stating our main result we will make a few remarks about cobordisms
between {\it non}-compact manifolds: 
let $(M,\omega,\phi)$ be a $2d$-dimensional Hamiltonian $S^1$-space
with moment map $\phi$.
As in the proof of Theorem 1 we will let
$S^1$ act on $M \times \CC$ diagonally, by the product of its action on $M$ 
and its standard action on $\CC$. If one equips $\CC$ with the two-form 
$\sqrt{-1} \, dz \wedge d \overline z$, 
this becomes a Hamiltonian action with moment map
$$
\psi (m, z) = \phi(m) + |z|^2\ .
$$
It is clear that $a$ is a regular value of $\phi$ if and only if $a$ is a
regular value of $\psi$, in which case one can reduce $M \times \CC$ at $a$ 
to obtain an orbifold,
$$
M^a : = \psi^{-1}(a)/S^1\ .
$$

\begin{lemma}
The orbifold $M^a$ is the disjoint union of the set
$$
M_{\phi < a} = \{ p \in M, \phi(p) < a\}
$$
and the reduced space $M_a = \phi^{-1}(a)/S^1$.
\end{lemma}

\begin{pf}
A pair $(m,z)$ is in $\psi^{-1}(a)$ if and only if $\phi(m) + |z|^2=a$. 
Let $z \not= 0$. Then $\phi(m) < a$ and
$$
z = e^{i\theta} (a - \phi(m)),\ \theta = \arg z\ .
$$
On the other hand, when $z=0$, we have $\phi(m)=a$, and the $S^1$-orbit 
through $(m, z)$ can be identified with a point of $M_a$. 
\hfill\end{pf}

In particular, if $\phi$ is proper and bounded from below, $M^a$ is
compact and hence is an orbifold compactification of the open subset
$M_{\phi<a}$ of $M$.  Notice also that, as above, the product of the
action of $S^1$ on $M$ with the trivial action on $\CC$ commutes with
the action we have just described and hence induces a Hamiltonian action
of $S^1$ on $M^a$. Moreover, if a compact Lie group
$G$ acts in a Hamiltonian fashion on $M$
and this action commutes with the action of $S^1$, one gets an 
induced Hamiltonian action of $G$ on $M^a$.

The operation
$$
M \mapsto M^a
$$
is called {\it symplectic cutting}. For more details about it see [Ler].

Now let $(M_r, \omega_r)$, for $r=1, 2$, be a Hamiltonian $S^1$-space with
moment map $\phi_r \colon M_r \to \RR$. 
We will assume that $\phi_r$ is proper and bounded from below. 

\begin{Definition}
The orbifolds $M_1$ and $M_2$ are cobordant as Hamiltonian $S^1$-spaces
if the cut spaces $M_1^a$ and $M_2^a$ are cobordant as Hamiltonian
$S^1$-spaces for all values of $a$.  
\end{Definition}

Suppose that, in addition, there is
a Hamiltonian action of $G$ on $M_r$ which
commutes with the action of $S^1$. Then we say that $M_1$ and $M_2$ are
cobordant as Hamiltonian $G$-spaces if $M_1^a$ and $M_2^a$
are cobordant as Hamiltonian $G$-spaces for all $a$.

\begin{Remarks}
\item
Although $M_r$ might be non-compact, the condition on the moment map
$\phi_r$ guarantees that the cut space $M_r^a$ is compact.
The cobording manifold between these cut spaces 
is required to be compact too.

\item
If $M_1^a$ is cobordant to $M_2^a$
for some value $a$ then  
$M_1^b$ is cobordant to $M_2^b$
for all $b<a$; this cobordism is obtained by taking the cobording
manifold with boundary $M_1^a \sqcup (-M_2^a)$, and cutting
it at the value $b$.

\item
One might hope to find a single (possibly non-compact) cobordism
$W$ between $M_1$ and $M_2$; its cuts, $W^a$, would then be
cobordisms between $M_1^a$ and $M_2^a$. However, this is not
always possible. 
\end{Remarks}

Let $(M, \omega, \phi)$ 
be a Hamiltonian $S^1$-space with moment map 
$\phi \colon M \to \RR$. We will assume that $\phi$ is proper and
bounded from below.  Also, for the moment we will assume that the 
fixed point set is 
discrete; $M^{S^1} = \{ p_1, p_2, \ldots \}$. 
For each fixed point $p_k$, there is an $\RR$-linear orientation
preserving map
\begin{equation} \labell{T-pk}
	T_{p_k} \cong \CC^d
\end{equation}
which converts the isotropy action of $S^1$ into the action
\begin{equation}\labell{8}
e^{i\theta}z = (e^{i m_{1k} \theta} z_1, \ldots, e^{i m_{dk}\theta} z_d ).
\end{equation}

If $\omega$ is symplectic, we can assume that the isomorphism
\eqref{8} converts the symplectic form on $T_{p_k}$ into the standard 
symplectic form on $\CC^d$:
\begin{equation}\labell{9}
	{\sqrt{-1} \over 2} \sum dz_r \wedge d \overline z_r\ .
\end{equation}
In any case, let $\omega_k$ be the symplectic form
\begin{equation}\labell{9a}
	\omega_k = {\sqrt{-1} \over 2}
	\sum \epsilon_{rk} dz_r \wedge d \overline z_r
\end{equation}
where $\epsilon_{rk} = \hbox{sgn}(m_{r,k})$, and let $\sigma_k$
be the integer
\begin{equation}\labell{10}
	\sigma_k = \# \{ m_{rk} \, , \, \epsilon_{rk} = -1 \}.
\end{equation}
The action \eqref{8} is Hamiltonian with respect to both the forms
above; however, the associated moment maps are different. The moment map
associated with the form \eqref{9} is 
$$
	{1 \over 2} \sum m_{rk} |z_k|^2 
	\quad \hbox{(plus an additive constant)},
$$
and the moment map associated with the form \eqref{9a} is
$$
	{1 \over 2} \sum \epsilon_{rk} m_{rk} |z_k|^2 
	\quad \hbox{(plus an additive constant).}
$$
Since $\epsilon_{rk} m_{rk} = |m_{rk}| >0$, the second of these maps is 
bounded from below and proper.

We can now state our main result.

\begin{theorem}[Linearization theorem]  \labell{linearization}
$M$ is cobordant, as a Hamiltonian $S^1$-space, 
to the disjoint union of the linear spaces $(\CC^d, \omega_k)$,
$k=1,2,\ldots$,
where $(\CC^d, \omega_k)$ is equipped with the Hamiltonian
action of $S^1$ defined by \eqref{8}, with the moment map
\begin{equation}\labell{12}
\phi_k = {1 \over 2} \sum \epsilon_{rk} m_{rk} |z_k|^2 + \phi(p_k)\ ,
\end{equation}
and with its complex orientation.
\end{theorem}

Note that the complex orientation on $\CC^d$ is equal to
$(-1)^{\sigma_k}$ times the symplectic orientation induced by 
$(\omega_k)^d$.

\begin{Remarks}
\item If, in addition, a compact Lie group $G$ acts on $M$ in a Hamiltonian
fashion, and if this action commutes with the action of $S^1$, then
this cobordism will be a cobordism of Hamiltonian $G$-spaces. 
(For instance, this will be the case if $S^1$ is a subgroup of 
the center of $G$.)
This shows that the circle in Theorem \ref{linearization} 
can be replaced by a torus acting on $M$ with isolated fixed points.

\item 
If the set of fixed points is not discrete, an analogue of 
Theorem \ref{linearization}
is true with the $T_{p_k}$'s replaced by the normal bundles to the 
fixed point components; for details see [GGK].

\item
If $\omega$ is not symplectic and $M$ is oriented,
the isotropy weights $m_{rk}$ for every fixed point $p_k$ 
are only determined up to a 
simultaneous change of sign of an even number of them.
However, the integers $\sigma_k$ and $\epsilon_{rk} m_{rk}$
are well defined and Theorem \ref{linearization} remains true
for any choice of $m_{rk}$'s.
\end{Remarks}

The rest of this article is devoted to applications of Theorem
\ref{linearization}. However, before getting into details, let us explain
the main idea of
how to use cobordisms to evaluate certain invariants of
Hamiltonian spaces. The cohomological invariants such as, for example, the
Duistermaat-Heckman measure
are, by Stokes's theorem, invariants
of cobordism.  This is also true for the cobordism class of the symplectic
reduction and for the (equivariant) Riemann-Roch number, i.e., the virtual
geometric quantization. (In the latter case, the notion of cobordism is
to be modified to take into account the stable complex structure of the
symplectic manifold; see Application 4 below.) By Theorem \ref{theorem1}
(or Theorem \ref{linearization}), a 
Hamiltonian space $M$ with isolated fixed points is cobordant
to a disjoint union of twisted projective spaces (or linear spaces)
associated with the local data near the fixed point set. Therefore, a
cobordism invariant of $M$ can be expressed as the sum of invariants
of these spaces. For example, when dealing with the Liouville measure,
this procedure leads immediately to the Duistermaat-Heckman formula
(Application 2).

\section*{Application 1 -- Symplectic Reduction}
Suppose that a torus $T$ acts in a Hamiltonian fashion 
on $(M,\omega)$ with isolated fixed points $\{ p_k\}$. 
Let $\phi \colon M \to {\frak t}^*$ be the moment map 
associated with this action and let
$\phi_k \colon \CC^d \to {\frak t}^*$ be the moment map associated with
the linear isotropy action \eqref{8} of $T$ on $T_{p_k} \cong \CC^d$
and with the two-form $\omega_k$ which is given by \eqref{9a}.
This map is unique
up to an additive constant and we will fix this constant by requiring that
$\phi_k(0) = \phi(p_k)$. 

Suppose now that $a \in {\frak t}^*$ is a regular value
of $\phi$ and of the $\phi_k$'s. Then, 
by Theorem \ref{theorem4}, the reduced space
\begin{equation}\labell{13}
	M_\red = \phi^{-1}(a)/T
\end{equation}
with its reduced two-form is cobordant to the disjoint union of
\begin{equation}\labell{Mk}
M_k = \phi_k^{-1}(a) / T\ ,\qquad k=1, 2,  \ldots .
\end{equation}
These spaces are compact symplectic toric orbifolds (see \cite{LT}).
This proves:

\begin{theorem} \labell{toric}
$M_\red$ is cobordant to a disjoint union of compact symplectic 
toric orbifolds.
\end{theorem}

This result was also proved by Shaun Martin \cite{Mar};
he expressed $M_\red$ as a ``tower" of twisted projective
spaces (i.e., a bundle over a bundle over $\ldots$ etc.,
where the fibers are twisted projective spaces).

\begin{Remark}
Theorem \ref{toric} can be modified to also cover the case
where $a$ is a regular value for $\phi$ but not a regular
value for the $\phi_k$'s. In this case, $M_\red$ is still
cobordant to a disjoint union of compact toric orbifold,
but the two-forms on these might be degenerate. For details,
see \cite{GGK}.
\end{Remark}

\section*{Application 2 -- The Duistermaat-Heckman theorem}
Let $T$ be a torus and let $(M,\omega,\phi)$ be a compact
Hamiltonian $T$-space with isolated fixed points.
Let $(\CC^d,\omega_k,\phi_k)$ be as in Theorem 
\ref{linearization}. 

Let $\mu$ be the measure on Borel subsets of $M$ associated 
with the top-form $\omega^d / d!$, and let $\mu_k$ be the measure 
on $\CC^d$ associated with $(\omega_k)^d/d!$. 
Their push-forwards,
$$
\nu = \phi_* \mu \quad \hbox{and} \quad \nu_k = (\phi_k)_* \mu_k,
$$
are the Duistermaat-Heckman measures associated with the 
Hamiltonian spaces $(M,\omega,\phi)$ and $(\CC^d, \omega_k, \phi_k)$. 
Theorem \ref{linearization}, 
coupled with the fact that the Duistermaat-Heckman 
measure is a cobordism invariant of Hamiltonian $T$-spaces, implies that
\begin{equation} \labell{15}
	\nu = \sum_k \nu_k\ .
\end{equation}
Note that $\mu$ and $\mu_k$ are {\em signed\/} measures, even in the
symplectic case; the top form $(\omega_k)^d/d!$ might be negative
with respect to the (complex) orientation of $\CC^d$.
If, instead, we integrate with respect to the {\em symplectic\/}
orientations, we get the more familiar formula of 
Guillemin-Lerman-Sternberg:
$$ \nu = \sum_k (-1)^{\sigma_k} |\nu_k|$$
which involves the positive push-forward measures $|\nu_k|$.

This identity was proved by other means in \cite{GLS} (for compact spaces)
and in \cite{PW} (for certain non-compact Hamiltonian spaces).

\section*{Application 3 -- The Jeffrey-Kirwan localization theorem}

Let us show how to use cobordisms to obtain the Jeffrey-Kirwan 
localization theorem 
in the abelian case.
By Application 1, the reduction $M_\red$ is
cobordant to a disjoint union of toric varieties $M_k$ associated
with the fixed points of the action. Therefore, for a cohomology
class $c\in H^*_T(M)$, the integral of the restriction of $c$ to $M_\red$ 
is equal to the sum of integrals of its restrictions to the $M_k$'s. 
This is in fact 
the Jeffrey-Kirwan localization theorem,
as long as we do not
care how to carry out the integration over $M_k$'s explicitly.

Let us get more specific.
For a regular value $a$ of the moment map $\phi$, there is a canonical map
(the Kirwan map)
$$
\kappa_a \colon H_T^*(M) \to H^*(M_\red)
$$
which maps the equivariant cohomology of $M$ (with coefficients in $\CC$)
surjectively onto the ordinary cohomology of the reduced space
 $M_\red = \phi^{-1}(a) /T$. 
By definition, this map is the composite of the restriction mapping
$$
	H_T^*(M) \to H_T^*(\phi^{-1}(a))
$$
and the Cartan isomorphism
$$
	H_T^*(\phi^{-1}(a)) \cong H^*(\phi^{-1}(a)/T)\ .
$$
Also, for each $p_k$, one has a canonical homomorphism,
$$
	\sigma_k \colon H_T^*(M) \to H_T^*(\{p_k\}) 
	{\buildrel \cong \over \to} H_T^*(\{0\}) \cong H_T^*(\CC^d)
$$
the first arrow being the restriction map. In the term on the right,
the action of $T$ on $T_{p_k} = \CC^d$ is the isotropy action given 
by \eqref{8}. In addition, there is a Kirwan map
$$
\kappa_k \colon H_T^*(\CC^d) \to H^*(M_k)
$$
where $M_k$ is the reduced space \eqref{Mk}, a toric variety. 
Applying Stokes's theorem to the
cobordism described in Theorem \ref{toric},
one can deduce, for $c \in H_T^*(M)$, the
following identity:
\begin{equation}\labell{16}
\int_{M_\red} \kappa(c) = \sum_k \int_{M_k} \kappa_k \circ
\sigma_k (c)\ .
\end{equation}

\begin{Remarks}
\item This is a topological form of the Jeffrey-Kirwan
localization theorem ([JK]). 
Their version of this theorem 
is valid for nonabelian groups.
However, Shaun Martin \cite{Mar} has recently given a purely 
topological proof that the abelian version of the localization theorem 
implies the non-abelian version.

\item 
An explicit recipe for evaluating the terms on the right 
of \eqref{16} is given in \cite{GS}.

\item
We have implicitly assumed that $a$ is a regular values
for the $\phi_k$'s. This assumption can be avoided \cite{GGK}.

\item
As before, if we equip the $M_k$'s with  their {\em symplectic\/} 
orientations, the summands of \eqref{16} need to be taken with the
coefficients $(-1)^{\sigma_k}$. 
\end{Remarks}

\section*{Application 4 -- Quantization}

A {\em stable complex structure\/} on a manifold $M$ is,
by definition, a complex structure on the bundle
$TM \oplus \RR^\ell$ for some $\ell$. 
Two such structures, one on the bundle $TM \oplus
\RR^{\ell_1}$ and one on the bundle $TM \oplus\RR^{\ell_2}$, are said to be
equivalent if there exists, for some choice of $m_1$ and $m_2$, an isomorphism
of complex vector bundles
$$
TM \oplus \RR^{\ell_1} \oplus \CC^{m_1} \cong TM \oplus \RR^{\ell_2} \oplus
\CC^{m_2}\ .
$$
For example, an almost complex structure on $M$ can always be viewed as 
a stable complex structure. Furthermore,
on a symplectic manifold $(M, \omega)$ there is a canonical
stable complex structure. Namely on $TM$
itself there exists a complex structure which is compatible with $\omega$ and
this complex structure is unique up to isomorphism. 
If $\omega$ is invariant with respect to the action of a compact group $G$,
the almost complex structure can also be chosen $G$-invariant.

Given two compact
oriented manifolds $M_r$, $r=1, 2$, each equipped with a stable complex
structure $J_r$, one says that $(M_1, J_1)$ and $(M_2, J_2)$ are {\it
cobordant} if there exists a compact oriented manifold-with-boundary $W$ and
a stable complex structure $J$ on $W$ such that
\begin{equation}\labell{17}
\partial W = M_1 \cup (-M_2)
\end{equation}
and
\begin{equation}\labell{18}
\iota_r^* J \sim J_r,
\end{equation}
$\iota_r \colon M_r \to \partial W$ being the inclusion mapping. (The
equivalence \eqref{18} makes sense in view of the fact that $\iota_r^*(TM
\oplus \RR^\ell) = TM_r \oplus \RR^{\ell+1}$.) With these definitions one has
the following addendum to Theorem \ref{linearization}.

\begin{theorem} \labell{stable-complex}
Let us equip $TM$ with 
an invariant complex structure $J$ 
and choose isomorphisms \eqref{T-pk} which respect
this structure.
Then the cobordism described in Theorem \ref{linearization} is a
cobordism of stable complex structures, i.e., we have
\begin{equation} \labell{cobordism}
	(M,\omega,\Phi,J) \sim 
	\sum_k (\CC^d, \omega_k, \phi_k, i) 
\enspace ,
\end{equation}
where $i$ denotes the intrinsic complex structure on $\CC^d$.
\end{theorem}

\begin{Remark}
As before, if on $(\CC^d,\omega_k)$ we take the symplectic
orientation instead of the complex orientation, we 
must put the coefficient $(-1)^{\sigma_k}$ in front of it.
\end{Remark}

Recall that the cobordism between $M$ and the $\CC^d$'s described 
in Theorems \ref{linearization} and \ref{stable-complex} is, 
strictly speaking, a cobordism between spaces which are obtained from 
these by the ``symplectic cutting'' operation.
Therefore, one item which remains to be explained in the statement of 
Theorem \ref{stable-complex}
is how the stable complex structures on $M$ and on the $\CC^d$'s give rise
to stable complex structures on the spaces obtained from them by symplectic
cutting. Since symplectic cutting is just a special case of symplectic
reduction, the answer to this is provided by the following:

\begin{theorem} \labell{theorem8}
Let $T$ be a torus and $M$ a Hamiltonian $T$-space. Fix
a regular value $a$ of the
moment map $\phi \colon M \to {\frak t}^*$. Then a $T$-invariant stable
complex structure on $M$ induces a stable complex structure on the reduced
space $M_\red = \phi^{-1}(a)/T$.
\end{theorem}

\begin{pf}
Letting $Z = \phi^{-1}(a)$ and letting $\iota \colon Z \to M$ and 
$\pi \colon Z \to M_\red$ 
be the inclusion and projection maps, it is easy to see that
\begin{equation}\labell{19}
\iota^* TM = \pi^* TM_\red \oplus ({\frak t} \oplus {\frak t}^*)\ .
\end{equation}
Thus a $T$-invariant complex structure on $TM \oplus \RR^\ell$ induces a
$T$-invariant complex-structure on $\pi^*TM_\red \oplus ({\frak t} \oplus
{\frak t}^*) \oplus \RR^\ell$. Since $T$ is abelian, the action of $T$ on
${\frak t} \oplus {\frak t}^*$ is trivial, so this is equivalent to a 
complex structure on 
$TM_\red \oplus ({\frak t} \oplus {\frak t}^*) \oplus \RR^\ell$. 
\hfill\end{pf}

Given a compact manifold $M$ with a stable complex structure $J$ and a
complex line bundle $\LL \to M$, one defines the {\em Riemann-Roch 
number\/} of $(M, J, \LL)$ to be the integral
\begin{equation}\labell{20}
\int_M \exp c(\LL) \hbox{Todd}(M,J) \ ,
\end{equation}
Todd$(M,J)$ being the Todd class of the complex vector bundle 
$TM \oplus \RR^\ell$
and $c(\LL)$ being the Chern class of $\LL$. One can also define an
equivariant version of \eqref{20} (cf.\ \cite{BGV}) 
and an orbifold version of \eqref{20} 
(the Kawasaki Riemann-Roch number of $(M, J, \LL)$, 
cf.\ \cite{Kaw}.)\footnote{
Warning: the formula for the Kawasaki Riemann-Roch number
is more complicated than \eqref{20}.
}

Suppose now that $(M,\omega)$ is a pre-quantizable Hamiltonian 
$T$-space with a pre-quantum line bundle $\LL$.
Assume for the moment that $\omega$ is 
symplectic.\footnote{
A detailed treatment of quantization in the case that the two-form
$\omega$ is degenerate will appear in \cite{CKT}.
}
If $a \in {\frak t}^*$ is an integer lattice point then
the reduced space \eqref{13} is pre-quantizable. 
Let $\LL_\red$ be its pre-quantum line bundle, 
$J_\red$ an almost complex structure on $M_\red$ which is
compatible with its symplectic structure, 
and $\RRoch(M_\red) = \RRoch(M_\red, J_\red, \LL_\red)$. 
Similarly, for each fixed point $p_k$, let $M_k$ be the 
reduced space \eqref{Mk} and
$\RRoch(M_k) = \RRoch(M_k, J_k, \LL_k)$ 
where $\LL_k$ is the reduced pre-quantum line bundle and $J_k$ is
the stable complex structure associated via Theorem \ref{theorem8}
with the intrinsic complex structure on $\CC^d$. 
(Usually this will {\it not} be compatible with 
the symplectic form on $M_k$!) 
Applying the Jeffrey-Kirwan theorem to the
equivariant version of $c(\LL)$ Todd$(M,J)$, one gets the following
``quantization'' identity:
\begin{equation}\labell{RRnumbers}
\RRoch(M_\red) = \sum_k \RRoch(M_k)
\end{equation}
expressing the Riemann Roch number of the reduced space in terms of 
Riemann-Roch numbers of toric orbifolds.

\begin{Remark}
A stable complex structure $J_k$ does {\em not\/}
determine an orientation. The orientation on $M_k$ which
we take in \eqref{RRnumbers} is induced by the complex orientation
of $\CC^d$; see \cite{CKT} for details.
If, instead, we take the symplectic orientation determined
by $\omega_k$, we must insert the coefficients
$(-1)^{\sigma_k}$ in front of the summands in \eqref{RRnumbers}.
\end{Remark}

Note that, in principle, $\RRoch(M_k)$ only depends on the linear isotropy
representation of $T$ on $T_{p_k}$, though in practice its evaluation 
relies on some rather deep theorems of Brion-Vergne \cite{BV}, Cappell-Shaneson 
\cite{CS} and Guillemin \cite{gu}. In [GGK] we will give
an explicit formula for it in terms of a partition function involving the
weights of this representation.

\ 

Stable complex structures can be incorporated in the notion of
cobordism of Hamiltonian $G$-spaces to make the geometric quantization,
i.e., the equivariant Riemann-Roch number, into an invariant of cobordism.
Thus consider the category of formal Hamiltonian
$G$-spaces, i.e., oriented $G$-manifolds (or orbifolds) equipped with
an equivariant closed 2-form $\omega+\phi$ (or just its equivariant
cohomology class) and a $G$-equivariant 
stable complex structure $J$. 
Clearly, an oriented symplectic $G$-manifold with a fixed moment map
can be viewed as a formal Hamiltonian $G$-space. An example of a different
nature is the $\CC^d$'s from Theorem \ref{stable-complex}, where the
symplectic form need not be compatible with the complex structure. Two
such manifolds or orbifolds are said to be {\em strongly\/} cobordant
if they are cobordant as Hamiltonian $G$-spaces and the cobordism
$W$ can be chosen 
to carry a $G$-equivariant stable complex structure which extends
the structures on its boundary.
Using symplectic cutting one can extend
this notion to non-compact spaces as well. Theorem 
\ref{stable-complex}
thus claims that $(M,\omega,\phi,J)$ is strongly cobordant to 
$\sum (\CC^d,\omega_k,\phi_k,i)$ where the linear spaces are given
their complex orientations. A stable complex $G$-manifold with a
$G$-pre-quantum line bundle $\LL$ can be made into a formal Hamiltonian
$G$-space by replacing $\LL$ by its $G$-equivariant first Chern class.
With this definition we see that {\em the equivariant Riemann-Roch
number \eqref{20} is an invariant of strong cobordism}.

\end{document}